\documentclass[10pt]{article}

\DeclareMathAlphabet{\mathpzc}{OT1}{pzc}{m}{it}

\usepackage{amsmath}
\usepackage{amssymb}
\usepackage{amsfonts}
\usepackage{amsthm}
\usepackage{mathrsfs}
\usepackage{ifsym}
\usepackage{fontenc}
\usepackage{textcomp}
\usepackage{tipa}
\usepackage{enumerate}
\usepackage[pdftex]{color, graphicx}
\numberwithin{equation}{section}
\usepackage{hyperref}

\begin{document}
 
\title{{\bf Explicit congruences for the class equations $H_{-7p}(X)$ and $H_{-28p}(X)$.}}         
\author{Patrick Morton}        
\date{July 4, 2022}          
\maketitle

\section{Introduction.}

Let $p$ be a prime and define
\begin{align*}
K_{7p}(x) = \begin{cases} H_{-28p}(x), \ &p \equiv 3 \ (\textrm{mod} \ 4);\\
H_{-7p}(x) H_{-28p}(x), \ &p \equiv 1 \ (\textrm{mod} \ 4),
\end{cases}
\end{align*}
where $H_{-d}(x)$ is the class equation for the quadratic discriminant $-d$.  Recall that $H_{-d}(x)$ is the minimal polynomial over $\mathbb{Q}$ of the $j$-invariants of elliptic curves with complex multiplication by the order $\mathcal{O} = \textsf{R}_{-d}$ of discriminant $-d$ in the field $\mathbb{Q}(\sqrt{-d})$. \medskip

In this paper I will give a computational proof of the following theorem, which gives the factorization of $K_{7p}(x)$ modulo $p$ for primes $p > 1217$.  This is the analogue of theorems in \cite[Thms. 1.1, 1.2]{mor5} and \cite[Thm. 2.1]{mor4} for the primes $l = 2, 3$ and $5$ concerning the factorization of the analogous polynomials $K_{lp}(x)$ (mod $p$).  This theorem is needed for the proofs of Conjectures 2 and 3 in \cite{mor6}, in the same way that the results of \cite{mor5} were needed for \cite{mor1} and \cite{mor2}.  \medskip

Let $\Phi_7(X,Y)$ denote the classical modular equation for level $7$ and $J_p(x)$ the polynomial
\begin{equation}
J_p(x) \equiv \sum_{k=0}^{n_p}{{2n_p + \varepsilon \atopwithdelims ( ) 2k + \varepsilon}{2n_p-2k \atopwithdelims ( ) n_p-k}(-432)^{n_p-k}(x-1728)^k} \ \ (\textrm{mod} \ p),
\label{eqn:1.1}
\end{equation}
where
$$n_p = [p/12], \ \ \varepsilon = \frac{1}{2} \left(1-\left(\frac{-1}{p}\right)\right).$$
(See \cite{mor0}, \cite{mor5}.) 

\newtheorem{thm}{Theorem}

\begin{thm} If $p$ is a prime satisfying
$$p \in \{211, 241, 263, 271, 283\} \cup \{p: p > 300\} - \mathcal{E}$$
where $\mathcal{E}$ is the exceptional set
\begin{align*}
\mathcal{E} & = \{307, 313, 317, 337, 353, 397, 401, 449, 521, 547, 569, 587, 593,\\
& \ \ \ \ \ 617, 641, 653, 773, 797, 857, 881, 953, 977, 1097, 1193, 1217\},
\end{align*}
then we have the following canonical factorization over $\mathbb{F}_p$:
\begin{align*}
K_{7p}(X) &\equiv H_{-7}(X)^{2\epsilon_{7}} H_{-28}(X)^{2\epsilon_{28}} q_{171}(X)^{4\epsilon_{171}} \prod_{d \in \mathfrak{L}}{H_{-d}(X)^{4\epsilon_d}}\\
& \ \times \prod_{d \in  \mathfrak{Q}}{H_{-d}(X)^{4\epsilon_d}} \times \prod_{d \in \mathfrak{R} - \{171\}}{H_{-d}(X)^{4\epsilon_d}} \times \prod_{i}{(X^2+a_iX+b_i)^2},
\end{align*}
where
\begin{align*}
\mathfrak{L} & = \{3, 12, 19, 24, 27\},\\
\mathfrak{Q} & = \{20, 40, 48, 52, 75, 115, 187\},\\
\mathfrak{R} &= \{96, 132, 160, 171, 180, 192, 195\};
\end{align*}
and where the quadratics $X^2+a_iX+b_i$ in the final product are the irreducible quadratic factors of $\gcd(\Phi_7(x^p,x), J_p(x))$ (mod $p$) which are distinct from the quadratic factors in the other products.  Here $q_{171}(x)$ is a quadratic factor of $H_{-171}(x)$ given by equation (\ref{eqn:2.8}) (with $\sqrt{57}$ mod $p$ chosen so that $q_{171}(x)$ is irreducible), and the exponents $\epsilon_d = 0$ or $1$ are defined in (\ref{eqn:2.5}), (\ref{eqn:2.7}) and (\ref{eqn:2.12}) below.
\end{thm}

Each of the polynomials whose power occurs in this factorization is either irreducible (mod $p$) (for $d \in \mathfrak{L} \cup \mathfrak{Q} \cup \{7, 28, 171\}$) or a product of distinct irreducible quadratics (for $d \in \mathfrak{R} - \{171\}$) when $\epsilon_d = 1$, except for $H_{-24}(x)$, which is a product of two distinct linear factors.  By ``canonical'' it is meant that the different powers are all relatively prime to each other for the given primes.  Several examples are calculated in Section 3.  See the Appendix for the quadratic and quartic polynomials $H_{-d}(x)$.  The linear class equations are given in the factorization (\ref{eqn:2.1}) below.  Also see \cite[pp. 237-238]{co}.  \medskip

\section{Proof of Theorem 1.}

The modular equation $\Phi_7(X,Y)=0$ can be computed using the resultant
\begin{align*}
7^{14} \Phi_7(x,y) = \textrm{Res}_z&((z^2-3z+9)(z^2-11z+25)^3-x(z-8),\\
&(z^2-3z+9)(z^2+229z+505)^3-y(z-8)^7).
\end{align*}
See \cite[p. 17]{mor6}.  Calculating on Maple, we have
\begin{align}
\notag \Phi_7(x,x) & = -x^2(x+15^3)(x-255^3)(x-2 \cdot 30^3)^2(x+96^3)^2 (x+ 3 \cdot 160^3)^2\\
\notag & \ \ \times (x^2-4834944x+14670139392)^2\\
\label{eqn:2.1} & = -H_{-3}(x)^2 H_{-7}(x)H_{-28}(x) H_{-12}(x)^2 H_{-19}(x)^2 H_{-27}(x)^2 H_{-24}(x)^2
\end{align}
and
\begin{equation}
\textrm{disc}_y \Phi_7(x,y) = -7^7 H_{-4}(x)^4 \prod_{d \in \mathfrak{L} \cup \mathfrak{Q} \cup \mathfrak{R}}{H_{-d}(x)^{e_d}},
\label{eqn:2.2}
\end{equation}
where
\begin{align*}
\mathfrak{L} & = \{3, 12, 19, 24, 27\},\\
\mathfrak{Q} & = \{20, 40, 48, 52, 75, 115, 187\},\\
\mathfrak{R} &= \{96, 132, 160, 171, 180, 192, 195\};
\end{align*}
and where $e_{3} = 6, e_{4} = 4$ and all other $e_d = 2$ ($d \neq 3, 4$).  (See \cite[p. 338]{f}.)  The polynomials $H_{-d}(x)$, for $d \in \mathfrak{L} \cup \{7, 28\}$ are the class equations which contribute linear factors to the factorization of $K_{7p}(x)$.  The class equations with $d \in \mathfrak{Q}$ are quadratic and irreducible (mod $p$) when they occur; and the class equations with $d \in \mathfrak{R}$ are all quartic polynomials, which for $d \neq 171$ factor into products of irreducible quadratics (mod $p$).  \bigskip

From \cite[Lemma 2.3]{mor5} we know that for $p > 28$ the irreducible factors of $K_{7p}(x)$ are the same as the supersingular factors which divide $\Phi_7(x^p,x)$, and their multiplicities are the same as their multiplicities in $\Phi_7(x^p,x)^2$:
\begin{equation}
K_{7p}(x) = \prod_{i}{q_i(x)^{e_i}}, \ \ \textrm{over} \ \ q_i(x) \mid \textrm{gcd}(ss_p(x),\Phi_7(x^p,x))
\label{eqn:2.3}
\end{equation} 
and $q_i(x)^{e_i} || \Phi_7(x^p,x)^2$.  For the proof of the theorem we shall calculate $F(t) = \Phi_7(t^p,t)$ and its first and second derivatives mod $p$ evaluated at a root $t$ of each supersingular factor $H_{-d}(X)$, for $d \in \mathfrak{L} \cup \mathfrak{Q} \cup \mathfrak{R} - \{171\}$, to show that $H_{-d}(X)^4 \ || \ K_{7p}(X)$ mod $p$.  For $d = 171$ it turns out that only one of the quadratic factors of the quartic polynomial $H_{-171}(x)$ over $\mathbb{Q}(\sqrt{57})$ can divide $K_{7p}(x)$ (mod $p$), but it does so to the fourth power, when it occurs.  These are the only factors which can occur to a power greater than $2$ in $K_{7p}(X)$, because
$$\textrm{disc}_x(\Phi_7(x^p,x)) \mid \Delta = \textrm{disc}_y(\Phi_7(x,y)).$$
See \cite[Prop. 2.4]{mor5}.  The calculations are most easily done using the de-symmetrized form of $\Phi_7(X,Y)$; i.e., the polynomial $Q_7(u,v)$, for which
$$Q_7(-x-y,xy) = \Phi_7(x,y).$$
The polynomial $Q_7(u,v)$ is given in the Appendix. \bigskip

\noindent {\it Linear factors.} \bigskip

With $Q(u,v)=Q_7(u,v)$ let
$$Q_1 = \frac{\partial Q(u,v)}{\partial u}, \ \ Q_2 = \frac{\partial Q(u,v)}{\partial v},$$
denote the first partials of $Q$ and $Q_{ij}$ the second partials.  If
$$F(t) = \Phi_7(t^p,t)=Q(-t^p-t,t^{p+1}),$$
then in characteristic $p$,
\begin{align*}
F'(t) & = -Q_1(-t^p-t, t^{p+1}) + t^p Q_2(-t^p-t, t^{p+1}),\\
F''(t) & = Q_{11}(-t^p-t, t^{p+1}) - 2t^pQ_{12}(-t^p-t, t^{p+1}) + t^{2p}Q_{22}(-t^p-t, t^{p+1}).
\end{align*}
For the roots $t$ of the linear factors we have
\begin{align*}
F(t) & =Q(-2t,t^2), \ \ F'(t) =-Q_1(-2t,t^2)+tQ_2(-2t,t^2),\\
\notag F''(t) & = Q_{11}(-2t,t^2)-2tQ_{12}(-2t,t^2) + t^2Q_{22}(-2t,t^2).
\end{align*}
For the discriminants occurring in the top part of Table \ref{tab:1}, $F(t) = F'(t) = 0$ in characteristic zero, so that at least $H_{-d}(x)^4$ occurs in the factorization of $K_{7p}(x)$ (mod $p$).  On the other hand, for $d = 7, 28$ we have $F(t) = 0$; but $F'(t)$ is only divisible by the primes listed in Table \ref{tab:1}, so only $H_{-7}(x)^2, H_{-28}(x)^2$ occur in the factorization for $p > 173$. \medskip

Note also that $H_{-4}(x) = x-1728$ does not occur at all in the factorization for $p > 23$, since
$$F(1728) = 2^{72} \cdot 3^{52} \cdot 7^2 \cdot 11^8 \cdot 19^4 \cdot 23^8.$$
Furthermore, the linear class equations for $d \in \mathfrak{L} \cup \{7, 28\}$ are pairwise relatively prime (mod $p$) for the primes not appearing in Table \ref{tab:1}. \medskip

\begin{table}
  \centering 
  \caption{Linear factors: calculating $F''(t)$ or $F'(t) \not \equiv 0$ (mod $p$).}\label{tab:1}

\noindent \begin{tabular}{|c|cl|c|}
\hline
    &  & \\
$d$	&   \ \ $t$  &  $F''(t)$   \\
\hline
   &  &  \\
3  & 0 & $2^{91} \cdot 3^{27} \cdot 5^{18} \cdot 7^3 \cdot17^9$  \\
  &  &  \\
12  &  $54000$ & $2^{39}  \cdot 3^{31} \cdot 5^{18} \cdot 7 \cdot11^4 \cdot 17^5 \cdot 23^2 \cdot 41^2 \cdot 47^2 \cdot 53^2 \cdot 59$  \\
  &  &  $\ \ \ \ \ \ \cdot 71^2 \cdot 83$\\
  &  &  \\
19  & $-96^3$ & $-2^{97} \cdot 3^{39} \cdot 7 \cdot 13^6 \cdot 29^2 \cdot 31 \cdot 41^2 \cdot 53^2 \cdot 89^2 \cdot 97 \cdot 113^2$ \\
  &  &   \\
 27 & $-12288000$ & $-2^{97} \cdot 3^9 \cdot 5^{19} \cdot 7 \cdot 11^4 \cdot 17^5 \cdot 23^4 \cdot 29^2 \cdot 41^3 \cdot 47 \cdot 89$\\
  &  &  $\ \ \ \ \ \ \cdot 113^2 \cdot 137^2 \cdot173$ \\
    &  &  \\
 24 & $2417472 \pm 1707264\sqrt{2}$ & $\alpha = A \pm B\sqrt{2}$\\
 & & $A = 2^{61} \cdot 3^{31} \cdot 7 \cdot 13^6 \cdot 17^5 \cdot 19^3 \cdot 23^2 \cdot 41 \cdot 47 \cdot 113^2$\\
 & &  $\ \ \ \ \ \ \cdot 2220151 \cdot 145642189738483$\\
 &  &  $B = 2^{64} \cdot 3^{33} \cdot 5^2 \cdot 7^3 \cdot 13^6 \cdot 17^5 \cdot 19^3 \cdot 23^2 \cdot 41 \cdot 47$\\
 & & $\ \ \ \ \ \ \cdot 113^2 \cdot 19553 \cdot 132578386723$\\
 & &  $A^2-2B^2 = 2^{122} \cdot 3^{62} \cdot 7^2 \cdot 13^{12} \cdot 17^{14} \cdot 19^6 \cdot 23^6 \cdot 41^4$\\
 & & $\ \ \ \ \ \ \cdot 47^3 \cdot 71^2 \cdot 89^2 \cdot 113^4 \cdot 137^2 \cdot 167$\\
  &  &  \\
\hline
   &  &  \\
$d$	&   \ \ $t$  &  $F'(t)$ \\
\hline
   &  &  \\
7 &  $-15^3$  &  $3^{46} \cdot 5^{21} \cdot 7 \cdot 13^7 \cdot 17^6 \cdot 19^4 \cdot 31^2 \cdot 41^4 \cdot 47^2$\\
 &  &  \\
 28 & $255^3$ & $-3^{44} \cdot 5^{23} \cdot 7 \cdot 13^7 \cdot 17^{12} \cdot 19^2 \cdot 41^2 \cdot 47^2 \cdot 59^2$\\
 &  &  $\ \ \ \ \ \  \cdot 83^2 \cdot 89^2 \cdot 97^2 \cdot 167^2 \cdot 173^2$\\
    &  &  \\
\hline
\end{tabular}
\end{table}

From these calculations it follows that the product of the linear factors which appear in the factorization of $K_{7p}(X)$ (mod $p$), for $p > 173$, is
\begin{equation}
L_p(X) = H_{-7}(X)^{2\epsilon_{7}} H_{-28}(X)^{2\epsilon_{28}} \prod_{d \in \mathfrak{L}}{H_{-d}(X)^{4\epsilon_d}},
\label{eqn:2.4}
\end{equation}
with
\begin{align}
\label{eqn:2.5} \epsilon_d & = \frac{1}{2}\left(1-\left(\frac{-d}{p}\right)\right), \ \ d \in \{3, 12, 19, 27\} \cup \{7, 28\};\\
\notag \epsilon_{24} &= \frac{1}{4}\left(1-\left(\frac{-24}{p}\right)\right) \left(1+\left(\frac{2}{p}\right)\right).
\end{align}

\begin{table}
  \centering 
  \caption{Quadratic factors: calculating $F''(t) \not \equiv 0$ (mod $p$).}\label{tab:2}

\noindent \begin{tabular}{|c|l|c|}
\hline
    &  \\
$d$	&    $\textrm{gcd}(D_1, D_2)$   \\
\hline
    &  \\
  $20$ & $2^{61} \cdot 5^{12} \cdot 7 \cdot 11^6 \cdot 13^5 \cdot 17^4 \cdot 19^2 \cdot 31^2 \cdot 37 \cdot 53^2 \cdot 73 \cdot 97 \cdot 113^2 \cdot 137$\\
     &  \\
 $40$  &  $2^{56} \cdot 3^{33} \cdot 5^{11} \cdot 7 \cdot 17^4 \cdot 31^2 \cdot 73 \cdot 83 \cdot 107 \cdot 113 \cdot 137^2 \cdot 163 \cdot 227 \cdot 233 \cdot 257$ \\
    &    \\
$48$ & $2^{21} \cdot 3^{26} \cdot 5^{16} \cdot 7 \cdot 11^6 \cdot 17^4 \cdot 29 \cdot 53^2 \cdot 89^2 \cdot 113 \cdot 137 \cdot 173 \cdot 197 \cdot 233 \cdot 257 $\\
   & $\cdot 269 \cdot 293 \cdot 317$\\
     &    \\
$52$ & $2^{61} \cdot 3^{32} \cdot 5^{16} \cdot 7 \cdot 13^3 \cdot 37 \cdot 41^2 \cdot 73 \cdot 89^2 \cdot 137 \cdot 193 \cdot 197 \cdot 281 \cdot 293 \cdot 317$\\
    &    $\cdot 353$\\
    &     \\
$75$ & $2^{82} \cdot 3^{27} \cdot 5^4 \cdot 7 \cdot 11^6 \cdot 17^4 \cdot 23 \cdot 29^2 \cdot 53^2 \cdot 83 \cdot 107 \cdot 113^2 \cdot 137 \cdot 233 \cdot 257$\\
     &    $\cdot 317 \cdot 353$\\
     &     \\
 $115$ & $2^{82} \cdot 3^{32} \cdot 5^{12} \cdot 7 \cdot 11^6 \cdot 13^5 \cdot 19^3 \cdot 23 \cdot 47 \cdot 61 \cdot 73 \cdot 163 \cdot 167 \cdot 173 \cdot 193 \cdot 197$\\
      &  $\cdot 257 \cdot 353 \cdot 593$   \\
      &     \\
 187 & $2^{76} \cdot 3^{33} \cdot 5^{16} \cdot 7 \cdot 11 \cdot 13^6 \cdot 17^4 \cdot 19^2 \cdot 23 \cdot 37 \cdot 71 \cdot 83 \cdot 113^2 \cdot 163 \cdot 269 \cdot 317$\\
   &   $\cdot 401 \cdot 521 \cdot 617 \cdot 641 \cdot 653 \cdot 881 \cdot 1193$\\
       &     \\  
\hline
\end{tabular}
\end{table}

\noindent {\it Quadratic factors.} \bigskip

For the quadratic factors $x^2+ux+v$ of $\textrm{disc}_y(\Phi(x,y))$ and a root $t$ of this factor, we have $F(t) = Q(u,v)$ and the formulas
\begin{align*}
F'(t) & = -Q_1(u, v) + t^p Q_2(u,v),\\
F''(t) & = Q_{11}(u, v) - 2t^pQ_{12}(u, v) + t^{2p}Q_{22}(u, v);
\end{align*}
and hence
$F''(t) = D_1(u,v)-t^p D_2(u,v)$, where
$$D_1(u,v)= Q_{11}(u,v)-vQ_{22}(u,v), \ \ D_2(u,v) = 2Q_{12}(u,v)+uQ_{22}(u,v).$$

For all of the discriminants occurring in Table \ref{tab:2}, $H_{-d}(x) = x^2+ux+v$ satisfies $Q(u,v) = Q_1(u,v) = Q_2(u,v) = 0$ in characteristic zero, so that $H_{-d}(x)^4$ divides $K_{7p}(x)$ (mod $p$) whenever $H_{-d}(x)$ is irreducible and its roots are supersingular in characteristic $p$.  Since $1, t^p$ are linearly independent over $\mathbb{F}_p$, for a root $t$ of one of the quadratics $H_{-d}(x)$, for $d \in \mathfrak{Q} - \{24\}$, it suffices to show $p \nmid \textrm{gcd}(D_1(u,v), D_2(u,v))$ to show that $F''(t) \not \equiv 0$ (mod $p$).  \medskip

It follows that for $p > 1193$, and all the primes not appearing in Table \ref{tab:2}, the contribution to $K_{7p}(x)$ (mod $p$) of the quadratic factors of the discriminant in (\ref{eqn:2.2}) is
\begin{equation}
Q_p(X) = \prod_{d \in  \mathfrak{Q}}{H_{-d}(X)^{4\epsilon_d}},
\label{eqn:2.6}
\end{equation}
where
\begin{equation}
\epsilon_d = \frac{1}{4}\left(1-\left(\frac{-d}{p}\right)\right) \left(1-\left(\frac{\textrm{disc}(H_{-d}(x))}{p}\right)\right).
\label{eqn:2.7}
\end{equation}
The exponential factor $\epsilon_d$ is chosen so that the quadratic $H_{-d}(x)$ appears in the product exactly when it is irreducible (mod $p$) and supersingular in characteristic $p$.  In addition, the quadratic class equations appearing in this product are relatively prime to each other (mod $p$), for primes not appearing in Table \ref{tab:2}.  \medskip

Note also that $H_{-24}(x) = x^2+ux+v$ does not occur in the factorization of $K_{7p}(x)$ when it is irreducible (mod $p$), for $p > 163$, because
\begin{equation*}
Q(u,v) = 2^{78} \cdot 3^{24} \cdot 13^8 \cdot 17^3 \cdot 19^5 \cdot 23^4 \cdot 37^3 \cdot 41^2 \cdot 43^2 \cdot 47 \cdot 61 \cdot 67 \cdot 109 \cdot 139 \cdot 157 \cdot 163.
\end{equation*}

\noindent {\it Quartic factors.} \bigskip

The calculations for all but one of the quartic factors in (\ref{eqn:2.2}) are similar to the calculations for the class equations $H_{-84}(x)$ and $H_{-96}(x)$ in \cite[pp. 104-107]{mor4}.  For $d \in \mathfrak{R}-\{171\}$, the Galois group of $H_{-d}(x)$ is $\mathbb{Z}_2 \times \mathbb{Z}_2$, and the same type of arguments apply.  However, the Galois group of $H_{-171}(x)$ is $D_4$, which requires a special argument.  This fact points out the peculiarities in the arguments for one prime versus another, i.e., in the current case $K_{7p}(x)$ versus the case $K_{5p}(x)$ discussed in \cite{mor4}. \medskip

For the discriminants for which $\textrm{Gal}(H_{-d}(x)/\mathbb{Q}) = \mathbb{Z}_2 \times \mathbb{Z}_2$, there are two primes $q_1, q_2$ dividing $d$, for which $\sqrt{q_1}, \sqrt{q_2}$ generate the splitting field of $H_{-d}(x)$ over $\mathbb{Q}$.  For each combination of Legendre symbols
$$\left(\left(\frac{q_1}{p}\right),\left(\frac{q_2}{p}\right)\right) \neq (+1, +1),$$
$H_{-d}(x) = f_1(x) f_2(x)$ factors as a product of two irreducible quadratics over a quadratic field $\mathbb{Q}(\sqrt{\delta})$, which is either $\mathbb{Q}(\sqrt{q_1}), \mathbb{Q}(\sqrt{q_2})$, or $\mathbb{Q}(\sqrt{q_1 q_2})$.  For exactly one of these values of $\delta$, its corresponding quadratics $f(x) = x^2+ux+v$ satisfy $Q(u,v) = Q_1(u,v) = Q_2(u,v) = 0$ in characteristic zero.  Table \ref{tab:3} contains the factorization of $\Delta= (N(D_1(u,v)),N(D_2(u,v)))$ for these two conjugate quadratic factors, where $N$ denotes the appropriate norm, showing that $H_d(x)^4$ exactly divides $K_{7p}(x)$ (mod $p$) when $p \nmid \Delta$ and $f(x) \in \mathbb{F}_p[x]$.  The second column contains the value $\delta = q_1, q_2$ or $q_1q_2$ whose square-root occurs in the factor $f(x)$.\medskip

\begin{table}
  \centering 
  \caption{Quadratic factors of quartics: calculating $F''(t) \not \equiv 0$ (mod $p$).}\label{tab:3}

\noindent \begin{tabular}{|c|c|l|c|}
\hline
    &  & \\
$d$	&  $\delta$ &  $\textrm{gcd}(N_{\mathbb{Q}(\sqrt{\delta})/\mathbb{Q}}(D_1), N_{\mathbb{Q}(\sqrt{\delta})/\mathbb{Q}}(D_2))$   \\
\hline
    &   &  \\
  $96$ & $2$ & $2^{68} \cdot 3^{52}\cdot 7^2 \cdot 13^{10}\cdot 17^8 \cdot 19^4 \cdot 37^2 \cdot 41^3 \cdot 43^2 \cdot 61^2 \cdot 89^2 \cdot 109^2  \cdot 113^2$\\
  &  &  $\cdot 137^3 \cdot 233 \cdot 257^2 \cdot 353^2 \cdot 401 \cdot 449 \cdot 521 \cdot 569 \cdot 593 \cdot 617^2 \cdot 641$\\
     &   &  \\
$132$ & $11$ & $2^{108} \cdot 3^{52} \cdot 5^{31} \cdot 7^2 \cdot 11^6 \cdot 13^{10} \cdot 31^4 \cdot 53^3 \cdot 67^2 \cdot 89^3 \cdot 113^2 \cdot 137^3$\\
   &   &  $\cdot 257^2 \cdot 317 \cdot 353 \cdot 401 \cdot 449 \cdot 617 \cdot 641 \cdot 653 \cdot 773 \cdot 797 \cdot 881$\\
     &   &   \\
 $160$ &  $2$  &  $2^{65} \cdot 3^{64} \cdot 5^{24} \cdot 7^2 \cdot 17^{10} \cdot 29^4 \cdot 31^4 \cdot 43^2 \cdot 67^2 \cdot 73^2 \cdot 83^2 \cdot 97 \cdot 109^2$\\
   & &  $\cdot 113^3 \cdot 137^2 \cdot 193 \cdot 233^2 \cdot 257^2 \cdot 353^2 \cdot 617 \cdot 857 \cdot 953 \cdot 977^2 \cdot 1097$\\
     &   &   \\
 $180$ & $15$ & $2^{123} \cdot 5^{24} \cdot 7^2 \cdot 11^8 \cdot 13^{10} \cdot 17^{10} \cdot 19^4 \cdot 31^4 \cdot 37^2 \cdot 53^3 \cdot 113^2 \cdot 137^3 \cdot 197$\\
  &   &  $\cdot 233^2 \cdot 257^2 \cdot 293 \cdot 317^2 \cdot 353 \cdot 593^2 \cdot 617 \cdot 653 \cdot 773 \cdot 797 \cdot 953 \cdot 977$\\
   &   &  $\cdot 1097 \cdot 1193 \cdot 1217$\\
   &   &   \\
 $192$ & $2$ & $2^{21} \cdot 3^{52} \cdot 5^{30} \cdot 7^2 \cdot 11^{12} \cdot 17^{10} \cdot 23^2 \cdot 29^4 \cdot 41^3 \cdot 59^2 \cdot 89^3 \cdot 101^2 \cdot 113^3$\\
 &  &  $\cdot 137^3 \cdot 233 \cdot 257^2 \cdot 281 \cdot 353 \cdot 521 \cdot 569 \cdot 593 \cdot 617 \cdot 641^2 \cdot 857 \cdot 881$\\
 & &  $ \cdot 977 \cdot 1097^2 \cdot 1193 \cdot 1217$\\
   &   &   \\
 $195$ & $65$ & $2^{165} \cdot 3^{52} \cdot 5^{24} \cdot 7^2 \cdot 13^5 \cdot 19^4 \cdot 29^2 \cdot 31^4 \cdot 43^2 \cdot 47 \cdot 83 \cdot 101$\\
  & &  $\cdot 137^3 \cdot 167 \cdot 197^2 \cdot 227 \cdot 293^2 \cdot 317 \cdot 353^2 \cdot 593 \cdot 617^2 \cdot 773 \cdot 977 \cdot 1097$\\
  & & $\cdot 1217$\\
    &   &   \\
\hline
\end{tabular}
\end{table}

For the other values of $\delta$ in Tables \ref{tab:4} and \ref{tab:5}, the quadratic factors $f(x) = x^2+ux+v$ of $H_{-d}(x)$ over $\mathbb{Q}(\sqrt{\delta})$ satisfy $Q(u,v) \neq 0$.  The factorizations of $N(Q(u,v))$ and the gcd $(N(Q_1(u,v)), N(Q_2(u,v)))$ show that only the square of one or more of these quadratic factors divides $K_{7p}(x)$ for large enough primes $p$, and this can be absorbed into the final product of factors of the form $(x^2+a_i x+b_i)^2$ in Theorem 1.  In \cite{mor4} we called these the sporadic factors for the discriminant $-d$.  These only occur for finitely many primes $p$.  The sporadic primes $p > 300$ with $\left(\frac{\delta}{p}\right) = +1$ are listed in bold in Tables \ref{tab:4} and \ref{tab:5}.  These are the primes $p > 300$ which do not appear in the corresponding factorization in Table \ref{tab:6}. \bigskip

\begin{table}
  \centering 
  \caption{Sporadic quadratic factors of quartics.}\label{tab:4}

\noindent \begin{tabular}{|c|c|l|c|}
\hline
    &  & \\
$d$	&  $\delta$ &  $N_{\mathbb{Q}(\sqrt{\delta})/\mathbb{Q}}(Q(u,v))$   \\
\hline
    &   &  \\
  $96$ & $3$ & $2^{84} \cdot 3^{48} \cdot 13^{16} \cdot 17^{14} \cdot 19^{12} \cdot 23^4 \cdot 37^5 \cdot 41^2 \cdot 43^6 \cdot 47^7 \cdot 61 \cdot 67^4 \cdot 89^2$\\
   & &  $\cdot 109^2 \cdot 113^2 \cdot 157^3 \cdot 181 \cdot 191^2 \cdot 229 \cdot 277^2 \cdot 397^3 \cdot {\bf 421} \cdot {\bf 541} \cdot {\bf 613} \cdot {\bf 661}$\\
       &   &  \\
  &  $6$ & $2^{72} \cdot 3^{48} \cdot 13^{16} \cdot 17^{14} \cdot 19^9 \cdot 23^{10} \cdot 37^6 \cdot 41^2 \cdot 43^2 \cdot 61^6 \cdot 67^4 \cdot 71^2 \cdot 89^2$\\
  &   &  $\cdot 137^2 \cdot 139 \cdot 163^2 \cdot 211 \cdot 283 \cdot 307^2 \cdot {\bf 331}^2 \cdot {\bf 379} \cdot {\bf 499} \cdot {\bf 523} \cdot 547^2 \cdot {\bf 571}$\\
  &  &  $\cdot {\bf 619} \cdot {\bf 643}$\\
  &  &   \\
 $132$ & $3$ & $2^{172} \cdot 3^{48} \cdot 5^{48} \cdot 11^7 \cdot 13^{16}\cdot 31^8 \cdot 53^2 \cdot 61 \cdot 73^4 \cdot 83^2 \cdot 103^4 \cdot 131^2 \cdot 167^2$\\
 & &  $\cdot 193^4 \cdot 227 \cdot 239^2 \cdot 241 \cdot 277 \cdot {\bf 337} \cdot {\bf 349} \cdot {\bf 373}^2 \cdot {\bf 457} \cdot {\bf 541} \cdot {\bf 601}^2 \cdot {\bf 613}$\\
   &  &  $ \cdot {\bf 673} \cdot {\bf 733} \cdot {\bf 853} \cdot {\bf 877}$\\
   &  &   \\
   & $33$ & $2^{144} \cdot 3^{48} \cdot 5^{48} \cdot 11^{14} \cdot 13^{16} \cdot 31^4 \cdot 53^2 \cdot 61^4 \cdot 67^3 \cdot 73^4 \cdot 103^2 \cdot 107^2 \cdot 109^4$\\
   &  &  $\cdot 137^2 \cdot 163^3 \cdot 199^2 \cdot 223^2 \cdot 257^2 \cdot {\bf 331 \cdot 367 \cdot 463 \cdot 487^2 \cdot 499 \cdot 631 \cdot 643}$\\
   & &  ${\bf \cdot 691 \cdot 727 \cdot 751 \cdot 823 \cdot 883 \cdot 907}$\\
    &  &   \\
$160$ & $5$ & $2^{84} \cdot 3^{108} \cdot 5^{24} \cdot 17^6 \cdot 29^6 \cdot 31^{10} \cdot 43^6 \cdot 61^2 \cdot 73^2 \cdot 79^2 \cdot 83^6 \cdot 97^2 \cdot 101^3$\\
& &  $\cdot 107^4 \cdot 109^2 \cdot 149^2 \cdot 181 \cdot 229 \cdot 269^3 \cdot {\bf 461 \cdot 509^2 \cdot 701 \cdot 821^2 \cdot 941}$\\
& &  ${\bf  \cdot 1021 \cdot 1061 \cdot 1109}$\\
 &  &   \\
  &  $10$ &  $-2^{72} \cdot 3^{107} \cdot 5^{35} \cdot 17^6 \cdot 29^8 \cdot 31^3 \cdot 43^2 \cdot 61^4 \cdot 67^3 \cdot 71^2 \cdot 73^2 \cdot 83^2 \cdot 97^2 \cdot 101^4$\\
  & & $\cdot 107^3 \cdot 149^4 \cdot 163 \cdot 199 \cdot 227^3 \cdot {\bf 347 \cdot 443 \cdot 467 \cdot 563} \cdot 587^2 \cdot {\bf 643 \cdot 683^2}$\\
  & &  $\cdot {\bf 787 \cdot 827 \cdot 947}$\\
    &  &   \\
 $180$ & $5$ & $2^{144} \cdot 5^{24} \cdot 11^{20} \cdot 13^{16} \cdot 17^6 \cdot 19^{10} \cdot 31 \cdot 37^6 \cdot 53^2 \cdot 59^5 \cdot 71^2 \cdot 73^4 \cdot 79$\\
  & &  $\cdot 97^4 \cdot 113^2 \cdot 131^2  \cdot 137^2 \cdot 151^3 \cdot 157^4 \cdot 191^2 \cdot 199^3 \cdot  211 \cdot 251 \cdot 271^2 \cdot {\bf 331}$\\
  & &  $\cdot {\bf 379 \cdot 439 \cdot 499 \cdot 619 \cdot 691 \cdot 739 \cdot 751 \cdot 811 \cdot 859 \cdot 991 \cdot 1171 \cdot 1231}$\\
      &  &   \\
  & $3$ & $2^{172} \cdot 3^2 \cdot 5^{34} \cdot 11^{20} \cdot 13^{16} \cdot 17^6 \cdot 19^{10} \cdot 31^{12} \cdot 37^5 \cdot 53^2 \cdot 73^4 \cdot 79^4 \cdot 97^3$\\
  & &  $\cdot 113^2 \cdot 137^2 \cdot 139^4 \cdot 151^4 \cdot 157 \cdot 193^2 \cdot 277 \cdot 313^2 \cdot {\bf 373 \cdot 397 \cdot 433 \cdot 577}$\\
  & &  $ \cdot {\bf 613 \cdot 673 \cdot 757 \cdot 877 \cdot 997 \cdot 1033 \cdot 1093 \cdot 1153 \cdot 1213 \cdot 1237}$\\
      &  &   \\
  $192$ & $6$ & $2^{28} \cdot 3^{82} \cdot 5^{51} \cdot 11^{24} \cdot 17^6 \cdot 23^4 \cdot 29^6 \cdot 41^2 \cdot 47 \cdot 53^6 \cdot 59^6 \cdot 83^4 \cdot 131^4$\\
  & &  $\cdot 149 \cdot 173^2 \cdot 179^4 \cdot 197^3 \cdot 269^2 \cdot 293^2 \cdot 317^3 \cdot {\bf 461 \cdot 557 \cdot 653 \cdot 701 \cdot 773}$\\
  &  &  $\cdot {\bf 797 \cdot 821 \cdot 941^2 \cdot 1013 \cdot 1061 \cdot 1181 \cdot 1277 \cdot 1301}$\\
      &  &   \\
      &  $3$ & $2^{24} \cdot 3^{84} \cdot 5^{52} \cdot 11^{15} \cdot 17^6 \cdot 23^4 \cdot 29^8 \cdot 41^2 \cdot 47^4 \cdot 53^4 \cdot 59^3 \cdot 71^4 \cdot 83^4$\\
      &  &  $\cdot 101^6 \cdot 107^4 \cdot 131 \cdot 149^4 \cdot 173^4 \cdot 179 \cdot 227^2 \cdot 251 \cdot {\bf 347 \cdot 419 \cdot 467 \cdot 491}$\\
     &   &  $\cdot {\bf 587 \cdot 683 \cdot 947 \cdot 971 \cdot 1019 \cdot 1163 \cdot 1187 \cdot 1283 \cdot 1307}$\\
     &   &   \\ 
\hline
\end{tabular}
\end{table}

\begin{table}
  \centering 
  \caption{Sporadic quadratic factors of quartics (continued).}\label{tab:5}

\noindent \begin{tabular}{|c|c|l|c|}
\hline
    &  & \\
$d$	&  $\delta$ &  $N_{\mathbb{Q}(\sqrt{\delta})/\mathbb{Q}}(Q(u,v))$   \\
\hline
    &   &  \\
 $195$ & $5$ & $2^{252} \cdot 3^{48} \cdot 5^{24} \cdot 13^{10} \cdot 19^9 \cdot 29^{10} \cdot 31^4 \cdot 43^6 \cdot 47^4 \cdot 83^2 \cdot 103^4 \cdot 109 \cdot 127^4$\\
   &   &  $\cdot 151 \cdot 229^2 \cdot 241^2 \cdot 269 \cdot 271 \cdot {\bf 331 \cdot 409 \cdot 421 \cdot 661 \cdot 709 \cdot 769 \cdot 1009}$\\
   &   &  ${\bf \cdot 1129 \cdot 1201 \cdot 1321}$\\
   &   &   \\
   & $13$ & $-2^{252} \cdot 3^{48} \cdot 5^{34} \cdot 13^8 \cdot 19^{12} \cdot 29^2 \cdot 31^8 \cdot 43^2 \cdot 47^4 \cdot 83^2 \cdot 101^2 \cdot 103 \cdot 127$\\
   &  &  $\cdot 131 \cdot 151^4 \cdot 157^2 \cdot 179^2 \cdot 251^2 \cdot 277^2 \cdot 283 \cdot {\bf 311} \cdot 313^2 \cdot 337^2 \cdot {\bf 419 \cdot 433}$\\
   &  &   $\cdot {\bf 673 \cdot 937 \cdot 997 \cdot 1093 \cdot 1153 \cdot 1297}$\\
     &   &   \\
\hline
\end{tabular}
\end{table}

\begin{table}
  \centering 
  \caption{Calculating $F'(t) \not \equiv 0$ for sporadic quadratic factors of quartics.}\label{tab:6}

\noindent \begin{tabular}{|c|c|l|c|}
\hline
    &  & \\
$d$	&  $\delta$ &  $\gcd(N_{\mathbb{Q}(\sqrt{\delta})/\mathbb{Q}}(Q_1(u,v)), N_{\mathbb{Q}(\sqrt{\delta})/\mathbb{Q}}(Q_2(u,v)))$   \\
\hline
    &   &  \\
  $96$ & $3$ & $2^{70} \cdot 3^{36} \cdot 13^{12} \cdot 17^{10} \cdot 19^8 \cdot 23^2 \cdot 37^2 \cdot 43^4 \cdot 47^2 \cdot 67^2 \cdot 109 \cdot 157 \cdot 191 \cdot 397$\\
       &   &  \\
  &  $6$ & $2^{62} \cdot 3^{36} \cdot 13^{12} \cdot 17^{10} \cdot 19^3 \cdot 23^5 \cdot 37^4 \cdot 43 \cdot 61^4 \cdot 67 \cdot 71 \cdot 307 \cdot 547$\\
  &  &   \\
 $132$ & $3$ & $2^{132} \cdot 3^{36} \cdot 5^{36} \cdot 11^2 \cdot 13^{12} \cdot 31^6 \cdot 73^2 \cdot 103^2 \cdot 131 \cdot 167 \cdot 193^2 \cdot 239$\\
   &  &   \\
   & $33$ & $2^{116} \cdot 3^{36} \cdot 5^{36} \cdot 11^{10} \cdot 13^{12} \cdot 31 \cdot 61^2 \cdot 67 \cdot 73^2 \cdot 107 \cdot 109^2 \cdot 163^2 \cdot 199 \cdot 223$\\
    &  &   \\
$160$ & $5$ & $2^{72} \cdot 3^{80} \cdot 5^{21} \cdot 29^2 \cdot 31^7 \cdot 43^4 \cdot 79 \cdot 83^4 \cdot 101 \cdot 107^2 \cdot 109 \cdot 149 \cdot 269$\\
 &  &   \\
  &  $10$ &  $2^{62} \cdot 3^{75} \cdot 5^{28} \cdot 29^6 \cdot 43 \cdot 61^2 \cdot 67 \cdot 71 \cdot 101^2 \cdot 149^2 \cdot 227 \cdot 587$\\
    &  &   \\
 $180$ & $5$ & $2^{116} \cdot 5^{21} \cdot 11^{13} \cdot 13^{12} \cdot 19^6 \cdot 37^4 \cdot 59^2 \cdot 71 \cdot 73^2 \cdot 97^2 \cdot 157^2 \cdot 191 \cdot 199$\\
      &  &   \\
  & $3$ & $2^{136} \cdot 3 \cdot 5^{26} \cdot 11^{13} \cdot 13^{12} \cdot 19^6 \cdot 31^8 \cdot 37^2 \cdot 73^2 \cdot 79^2 \cdot 97 \cdot 139^2 \cdot 151^2 \cdot 193$\\
    &   &  $ \cdot 313$\\
      &  &   \\
  $192$ & $6$ & $2^{24} \cdot 3^{61} \cdot 5^{40} \cdot 11^{18} \cdot 23^2 \cdot 29^2 \cdot 53^3 \cdot 59^4 \cdot 83^2 \cdot 131^2 \cdot 173 \cdot 179^2 \cdot 197$\\
  & &  $\cdot 269 \cdot 293 \cdot 317$\\
      &  &   \\
      &  $3$ & $2^{22} \cdot 3^{64} \cdot 5^{38} \cdot 11^4 \cdot 23^2 \cdot 29^6 \cdot 47 \cdot 53^2 \cdot 71^2 \cdot 83^2 \cdot 101^4 \cdot 107^2 \cdot 149^2$\\
      &  &  $\cdot 173^2 \cdot 227$\\
     &   &   \\
  $195$ & $5$ & $2^{190} \cdot 3^{36} \cdot 5^{21} \cdot 13^8 \cdot 19^3 \cdot 29^6 \cdot 31 \cdot 43^4 \cdot 47^2 \cdot 103^2 \cdot 127^2 \cdot 229 $\\
   &  &  \\
   &  $13$ & $2^{190} \cdot 3^{36} \cdot 5^{26} \cdot 13^6 \cdot 19^8 \cdot 31^6 \cdot 43 \cdot 47^2 \cdot 101 \cdot 151^2 \cdot 179 \cdot 251 \cdot 277$\\
   &  &  $ \cdot 313 \cdot 337$\\
    &  &  \\
\hline
\end{tabular}
\end{table}

\noindent {\it The class equation $H_{-171}(x)$.} \bigskip

The Galois group of the polynomial $H_{-171}(x)$ is $D_4$, as stated before, and its splitting field is $\Omega_3$, the ring class field of conductor $f = 3$ over $\mathbb{Q}(\sqrt{-19})$.  Its discriminant is
\begin{align*}
\mathfrak{d} = \textrm{disc}(H_{-171}(x)) = & \ -2^{186} \cdot 3^{21} \cdot 13^{12} \cdot 19^6 \cdot 29^2 \cdot 31^6 \cdot 37^2 \cdot 41^2\\
& \ \ \ \ \cdot 59^2 \cdot 71^2 \cdot 79^2 \cdot 103^2 \cdot 127^2 \cdot 151^2.
\end{align*}
By the Pellet-Stickelberger-Voronoi theorem (PSV, see \cite[Appendix]{brm}), $H_{-171}(x)$ has an even or an odd number of irreducible factors (mod $p$), according as $\left(\frac{\mathfrak{d}}{p}\right) = \left(\frac{-3}{p}\right) = +1$ or $-1$. It splits into linear factors (mod $p$) if and only if $p$ splits completely in $\Omega_3$, and this holds if and only if $\left(\frac{-19}{p}\right) = +1$ and the prime factors of $p = \lambda \lambda'$ in $\mathbb{Q}(\sqrt{-19})$ have the form
$$\lambda  = \frac{x+y\sqrt{-19}}{2}, \ \ \textrm{with} \ \ 3 \mid y.$$
Furthermore, when $\left(\frac{-19}{p}\right) = +1$, the degrees of the prime divisors $\mathfrak{p}$ of $p$ in $\Omega_f$ are:
\begin{align*}
\textrm{deg} \ \mathfrak{p} = 4 &, \ \textrm{if} \ (xy,3) = 1;\\
\textrm{deg} \ \mathfrak{p} = 2 &,  \ \textrm{if} \ 3 \mid x.
\end{align*}
On the other hand, when $\left(\frac{-19}{p}\right) = -1$, $p$ is inert in $\mathbb{Q}(\sqrt{-19})$ and splits completely in $\Omega_3/\mathbb{Q}(\sqrt{-19})$.  In this case the factors of $H_{-171}(x)$ are supersingular in characteristic $p$. \medskip

$H_{-171}(x)$ factors over the quadratic field $k = \mathbb{Q}(\sqrt{57})$.  Its factors over this field are
\begin{align}
\label{eqn:2.8} q_{171}(x) = & \ x^2 + (347141028938268672 - 45979952529358848\sqrt{57})x\\
\notag & \ \ - 53697050710005252096 + 7112348114438062080\sqrt{57}\\
\notag = & \ x^2 +ux+v
\end{align}
and its conjugate $\tilde q_{171}(x)$.  We have $Q(u,v) = Q_1(u,v) = Q_2(u,v) = 0$ in characteristic $0$ and 
\begin{align*}
\textrm{gcd}&(N_{k/\mathbb{Q}}(D_1(u,v)),N_{k/\mathbb{Q}}(D_1(u,v))) = 2^{158} \cdot 3^{27} \cdot 7^2 \cdot 13^{10} \cdot 19^2 \cdot 29 \cdot 31^4\\
& \cdot 37^2 \cdot 41^2 \cdot 53^2 \cdot 59 \cdot 67^2 \cdot 71 \cdot 89^2 \cdot 107 \cdot 109^2 \cdot 113^2 \cdot 167 \cdot 173 \cdot 227\\
& \cdot 257 \cdot 269 \cdot 281 \cdot 293 \cdot 317 \cdot 449 \cdot 569 \cdot 641 \cdot 797 \cdot 857 \cdot 953 \cdot 977 \cdot 1097.
\end{align*}
Hence, the factor $q_{171}(x)^4$ divides $K_{7p}(x)$ (mod $p$) whenever it is supersingular, irreducible and
\begin{equation}
\left(\frac{-3}{p}\right) = \left(\frac{-19}{p}\right) = -1,
\label{eqn:2.9}
\end{equation}
for the primes not dividing the above gcd.  By the PSV theorem quoted above, $H_{-171}(x)$ has an odd number of irreducible factors (mod $p$), so it splits as a product of a quadratic and two linear factors.  \medskip

To determine when the factor $q_{171}(x)$ is irreducible (mod $p$), we calculate
\begin{align*}
\textrm{disc}(q_{171}(x)) & = 241013787944639410404629344363216896\\
& \ \ \ \ - 31923056063148787797313080401068032\sqrt{57}\\
& = 2^{30} \cdot 3^4 \cdot 13^2 \cdot 19 \cdot \alpha\\
 \alpha & = \frac{1726023441844820110263 - 228617389777623914371\sqrt{57}}{2},
\end{align*}
where
$$N_{k/\mathbb{Q}}(\alpha) = -2^2 \cdot 3 \cdot 29^2 \cdot 41^2 \cdot 59^2 \cdot 71^2.$$
As above,
$$\left(\frac{N_{k/\mathbb{Q}}(\textrm{disc}(q_{171}(x)))}{p}\right) = \left(\frac{-3}{p}\right) = -1.$$
Hence, one of $q_{171}(x)$ and its conjugate $\tilde q_{171}(x)$ is irreducible (mod $p$) and one is reducible.  Therefore, when (\ref{eqn:2.9}) holds, only the fourth power of the (irreducible) factor $q_{171}(x)$ divides $K_{7p}(x)$ (mod $p$) instead of $H_{-171}(x)^4$.  Note that
\begin{align*}
\alpha & = -(151+20\sqrt{57})^{-6}\frac{(128562177 + 17631211\sqrt{57})}{2}\\
& = -(151+20\sqrt{57})^{-6} \left( \frac{7+\sqrt{57}}{2}\right)^2 (15+2\sqrt{57})(22-3\sqrt{57})^2\\
& \ \ \times (4-\sqrt{57})^2(13-2\sqrt{57})^2(29+4\sqrt{57})^2,
\end{align*}
where $\varepsilon = 151+20\sqrt{57}$ is the fundamental unit in $k$ and the numbers inside each of the other parentheses are primes in $k$.  (The field $k$ has class number $1$.)  Hence, $q_{171}(x)$ is irreducible mod $p$ if and only if
$$\left(\frac{-19}{p}\right) \left(\frac{15+2\sqrt{57}}{p}\right) = -1;$$
when the factor $q_1(x)$ is supersingular, this happens if and only if
$$\left(\frac{15+2\sqrt{57}}{p}\right) =+1,$$ 
for the choice of the square-root $\sqrt{57}$ modulo $p$.
\bigskip

For the remainder of the calculation, we need a root
\begin{align*}
\eta_1 & = -173570514469134336 + 22989976264679424\sqrt{57}\\
& \ \  + 958464\sqrt{378534796\sqrt{-3} + 32794445395051582094997}\\
& \ \  - 958464\sqrt{-378534796\sqrt{-3} + 32794445395051582094997}\\
& = -173570514469134336 + 22989976264679424\sqrt{57}\\
& \ \ + 958464 (-36178080131 + 31806062964\sqrt{-3})\sqrt{4\sqrt{-3} - 3}\\
& \ \  - 958464 (-36178080131 - 31806062964\sqrt{-3})\sqrt{-4\sqrt{-3} - 3}
\end{align*}
of $q_{171}(x) = (x-\eta_1)(x-\eta_2)$ and a root
\begin{align*}
\eta_3 & = -173570514469134336 - 22989976264679424\sqrt{57}\\
& \ \ + 958464 (-36178080131 + 31806062964\sqrt{-3}) \sqrt{4\sqrt{-3} - 3}\\
& \ \  + 958464(-36178080131 - 31806062964\sqrt{-3}) \sqrt{-4\sqrt{-3} - 3}
\end{align*}
of its conjugate $\tilde q_{171}(x) = (x-\eta_3)(x-\eta_4)$.  If $(x-\eta_1)(x-\eta_3) = x^2+ux+v$ and $L=\mathbb{Q}(\eta_1,\eta_3)$, then
\begin{align*}
N_{L/\mathbb{Q}}&(Q(u,v)) = 2^{508} \cdot 3^{48} \cdot 13^{32} \cdot 19^{11} \cdot 29^8 \cdot 31^{13} \cdot 37^{11} \cdot 41^6 \cdot 53^6 \cdot 59^4\\
& \cdot 67^3 \cdot 79^5 \cdot 89^4 \cdot 97^6 \cdot 103^6 \cdot 107^2 \cdot 109^2 \cdot 113^2 \cdot 127^4 \cdot 151^5 \cdot 173^2 \cdot 181^2\\
& \cdot 193^4 \cdot 211 \cdot 223^2 \cdot 241^2 \cdot 337^2 \cdot 409 \cdot 433 \cdot 601 \cdot 673 \cdot 829 \cdot 877 \cdot 1009\\
& \cdot 1021 \cdot 1117 \cdot 1129 \cdot 1153.
\end{align*}
Moreover,
\begin{align}
\notag \textrm{gcd}&(N_{L/\mathbb{Q}}(Q(u,v)),N_{L/\mathbb{Q}}(Q_1(u,v)) = 2^{442} \cdot 3^{44}\cdot 13^{24} \cdot 19^7 \cdot 29^6 \cdot 31^8 \cdot 37^6 \cdot 41^2\\
\label{eqn:2.10} & \cdot 53^2 \cdot 59^2 \cdot 67 \cdot 79^2 \cdot 97^2 \cdot 103^2 \cdot 109 \cdot 127^2 \cdot 151^2 \cdot 181 \cdot 193^2 \cdot 337.
\end{align}
(Note that $\left(\frac{57}{p}\right) = -1$ for the primes $p \ge 181$ in the first factorization.)  It follows that the quadratic $x^2+ux+v$ divides $K_{7p}(x)$ at most for the primes appearing in the first factorization, and the second shows that when it does divide, it does so to only the second power, for $p > 337$.  A similar calculation for $(x-\eta_1)(x-\eta_4) = x^2+ux+v$ leads to the same result.  Since the Galois group of $H_{-171}(x)$ is $D_4$, there is an automorphism of $\Omega_3/k$ taking $\eta_3$ to $\eta_4$ and fixing $\eta_1$ and $\eta_2$.  Thus, the norms to $\mathbb{Q}$ of the above quantities have to be the same when $\eta_3$ is replaced by $\eta_4$.  The same results hold for the quadratics $(x-\eta_2)(x-\eta_3)$ and $(x-\eta_2)(x-\eta_4)$, since complex conjugation interchanges $\eta_1$ and $\eta_2$ and fixes the real roots $\eta_3$ and $\eta_4$.  Thus, these factors are sporadic factors when $\left(\frac{-3}{p}\right) = +1$, and can be absorbed into the final product in the factorization given in Theorem 1. \medskip

It is remarkable that the prime factors of these quantities are relatively small (less than $1217$(!)), especially since small errors in the calculation lead to some very large prime factors. \medskip

Computing the prime factorizations of the resultants $\textrm{Res}(H_{-d_1}(x), H_{-d_2}(x))$ for $d_1, d_2 \in \mathfrak{Q} \cup \mathfrak{R}$ shows that the polynomials $H_{-d_1}(x), H_{-d_2}(x)$, which are quadratic or have irreducible quadratic factors with $\epsilon_{d_i} = 1$, have no common factors modulo $p$, if $p \in \{211, 241, 263, 271, 283\}$ or $p >300$ and $p \notin \mathcal{E}$. \medskip

This discussion shows that the contribution to the factorization of $K_{7p}(x)$ by the class equations $H_{-d}(x)$ for $d \in \mathfrak{R}$ and $p > 1217$ is exactly
\begin{equation}
R_p(X) = q_{171}(X)^{4\epsilon_{171}} \prod_{d \in \mathfrak{R} - \{171\}}{H_{-d}(X)^{4\epsilon_d}},
\label{eqn:2.11}
\end{equation}
where $q_{171}(x)$ is given by (\ref{eqn:2.8}) and the choice of $\sqrt{57}$ in (\ref{eqn:2.8}) and the exponents $\epsilon_d$ satisfy
\begin{align}
\notag \epsilon_{171} &= \frac{1}{8} \left(1-\left(\frac{-3}{p}\right)\right)  \left(1-\left(\frac{-19}{p}\right)\right) \left(1+\left(\frac{15+2\sqrt{57}}{p}\right)\right),\\
\label{eqn:2.12} \epsilon_d & = \frac{1}{8} \left(1-\left(\frac{-d}{p}\right)\right)\left(1+\left(\frac{2}{p}\right)\right) \left(1-\left(\frac{d/32}{p}\right)\right), \ d \in \{96, 160, 192\},\\
\notag \epsilon_d & = \frac{1}{8} \left(1-\left(\frac{-d}{p}\right)\right) \left(1+\left(\frac{\delta}{p}\right)\right) \left(1-\left(\frac{q}{p}\right)\right), \ q = 3, 5, 5, \ \textrm{for} \ d = 132, 180, 195.
\end{align}

Hence, by (\ref{eqn:2.4}), (\ref{eqn:2.6}) and (\ref{eqn:2.11}) we have
\begin{align*}
K_{7p}(X) &\equiv L_p(X) Q_p(X) R_p(X) \ (\textrm{mod} \ p),\\
& \equiv H_{-7}(X)^{2\epsilon_{7}} H_{-28}(X)^{2\epsilon_{28}} q_{171}(X)^{4\epsilon_{171}} \prod_{d \in \mathfrak{L}}{H_{-d}(X)^{4\epsilon_d}}\\
& \ \times \prod_{d \in  \mathfrak{Q}}{H_{-d}(X)^{4\epsilon_d}} \times \prod_{d \in \mathfrak{R} - \{171\}}{H_{-d}(X)^{4\epsilon_d}} \times \prod_{i}{(X^2+a_iX+b_i)^2},
\end{align*}
where the quadratics $X^2+a_iX+b_i$ in the final product are irreducible (mod $p$) and distinct from the quadratic factors in the other products.  They divide $\gcd(\Phi_7(x^p,x),J_p(x))$ by (\ref{eqn:2.3}), since $ss_p(x) = x^r(x-1728)^s J_p(x)$ is a product of linear factors times $J_p(x)$ (mod $p$).  This factorization (mod $p$) holds for all primes which do not appear in Tables \ref{tab:1}-\ref{tab:3} and \ref{tab:6} or in the factorizations just before (\ref{eqn:2.9}) and in (\ref{eqn:2.10}).  These are the primes $p \in \{211, 241, 263, 271, 283\} \cup \{p: p > 300\}$ which are not in the set
\begin{align*}
\mathcal{E} & = \{307, 313, 317, 337, 353, 397, 401, 449, 521, 547, 569, 587, 593,\\
& \ \ \ \ \ 617, 641, 653, 773, 797, 857, 881, 953, 977, 1097, 1193, 1217\}.
\end{align*}

This completes the proof of Theorem 1.  This factorization is similar to the factorization of $K_{5p}(X)$ (mod $p$) in \cite{mor4}, but differs in that only a partial factor of $H_{-171}(X)$ occurs in the factorization.  As we have seen, this is a consequence of the fact that the Galois group of this class equation is $D_4$ rather than $\mathbb{Z} \times \mathbb{Z}$.

\section{Examples.}

{\bf 1}. $p = 211$.  We have
\begin{align*}
J_{211}(x) &= (x + 13)(x + 63)(x + 97)(x + 129)(x + 183)(x^2 + 56x + 23)\\
& \ \ \times (x^2 + 152x + 88)(x^2 + 162x + 146)(x^2 + 121x + 206)\\
& \ \ \times (x^2 + 186x + 97)(x^2 + 184x + 86) \ \ (\textrm{mod} \ 211),
\end{align*}
and
$$\gcd(\Phi_7(x^{211},x), J_{211}(x)) = (x+13)(x^2 + 56x + 23)(x^2 + 152x + 88)(x^2 + 186x + 97)$$
modulo $211$.  In addition, all exponential factors $\epsilon_d = 0$ except for $\epsilon_{19}=1$, and $H_{-19}(x) = x+96^3 \equiv x+13$ (mod $211$).  Hence, Theorem 1 gives the congruence
$$H_{-28 \cdot 211}(x) \equiv (x+13)^4 (x^2 + 56x + 23)^2(x^2 + 152x + 88)^2(x^2 + 186x + 97)^2 \ (\textrm{mod} \ 211),$$
which agrees with $h(-28 \cdot 211) = 16$.  Note that the factor $x^2 + 152x + 88$ is a sporadic factor of $H_{-171}(x)$. \medskip

\noindent {\bf 2}. $p=241$.  In this case $\epsilon_7 = \epsilon_{28} = \epsilon_{19} = \epsilon_{52} = 1$ and all other $\epsilon_d = 0$.  We also have
\begin{align*}
J_{241}(x) &= (x + 1)(x + 25)(x + 148)(x + 177)(x + 213)(x + 233)(x^2 + 11x + 156)\\
& \ \ \times (x^2 + 27x + 44)(x^2 + 28x + 79)(x^2 + 111x + 25)(x^2 + 160x + 117)\\
& \ \ \times (x^2 + 166x + 180)(x^2 + 206x + 211) \ \ (\textrm{mod} \ 241)
\end{align*}
and
\begin{align*}
\gcd(\Phi_7(x^{241},x), &J_{241}(x)) = (x + 1)(x + 25)(x + 148)(x^2 + 11x + 156)\\
& \times (x^2 + 27x + 44)(x^2 + 28x + 79)(x^2 + 111x + 25)\\
& \times (x^2 + 160x + 117)(x^2 + 166x + 180)
\end{align*}
modulo $p = 241$.  This gives the congruence
\begin{align*}
H_{-7 \cdot 241}(x) & H_{-28 \cdot 241}(x) \equiv (x + 1)^2(x + 148)^2(x + 25)^4(x^2 + 160x + 117)^4\\
& \ \ \times (x^2 + 11x + 156)^2 (x^2 + 27x + 44)^2(x^2 + 28x + 79)^2\\
& \ \ \times (x^2 + 111x + 25)^2 (x^2 + 166x + 180)^2  \ \ (\textrm{mod} \ 241),
\end{align*}
agreeing with the fact that $h(-7 \cdot 241) = h(-28 \cdot 241) =18$.  In this case
$$H_{-19}(x) = x+96^3 \equiv x+25, \ \ \ H_{-52}(x) \equiv x^2 + 160x + 117 \ \ (\textrm{mod} \ 241).$$

\noindent {\bf 3}. $p=311$.  In this case $\epsilon_d = 1$ for
$$d \in \{3, 7, 28, 12, 27, 24, 187\}.$$
Further, we have
\begin{align*}
J_{311}(x) & \equiv (x + 7)(x + 12)(x + 59)(x + 64)(x + 79)(x + 86)(x + 114)(x + 143)\\
& \times (x + 161)(x + 179)(x + 180)(x + 209)(x + 212)(x + 234)(x + 265)\\
& \times (x + 279)(x + 292)(x^2 + 158x + 119)(x^2 + 198x + 61)(x^2 + 22x + 260)\\
& \times (x^2 + 148x + 243) \ \ (\textrm{mod} \ 311);
\end{align*}
while
\begin{align*}
\gcd(\Phi_7(x^{311},&x), J_{311}(x)) \equiv (x + 12)(x + 79)(x + 114)(x + 161)(x + 212)(x + 265)\\
& \ \ \times (x^2 + 22x + 260)(x^2 + 148x + 243)(x^2 + 158x + 119) \ \ (\textrm{mod} \ 311).
\end{align*}
Hence,
\begin{align*}
H_{7 \cdot 311}(x) &\equiv H_{-7}(x)^2 H_{-28}(x)^2 H_{-3}(x)^4 H_{-12}(x)^4 H_{-27}(x)^4\\
& \ \times H_{-24}(x)^4 H_{-187}(x)^4 (x^2 + 148x + 243)^2 (x^2 + 158x + 119)^2 \\
& \equiv (x+265)^2 (x+212)^2 x^4 (x+114)^4 (x+79)^4 (x+12)^4 (x+ 161)^4\\
& \ \ \times (x^2+22x+260)^4(x^2 + 148x + 243)^2 (x^2 + 158x + 119)^2 \ \ (\textrm{mod} \ 311).
\end{align*}
This agrees with the fact that $h(-7 \cdot 311) = 40$. \medskip

When primes lie in the set $\{p: p < 300\} \cup \mathcal{E}$, a congruence for $H_{-28p}(x)$ or $H_{-7p}(x) H_{-28p}(x)$
can be computed by finding the gcd of $\Phi_7(x^p,x)$ and $J_p(x)^e$ (mod $p$) for high enough exponents $e$.  We illustrate this in the next two examples.  Note that $H_{-3}(x)^4 = x^4$ exactly divides $K_{7p}(x)$ mod $p$ for $p > 17$ and $p \equiv 2$ (mod $3$), by the factorization in Table \ref{tab:1}. \medskip

\noindent {\bf 4}. $p = 113$.  Here we have
$$J_{113}(x) \equiv (x + 14)(x + 41)(x + 59)(x^2 + 9x + 82)(x^2 + 18x + 38)(x^2 + 32x + 65)$$
modulo $113$, and
\begin{align*}
\Phi_7(x^{113},x) &\equiv 55(x + 14)^2 (x + 41)^4 (x + 59)^4 (x^2 + 9x + 82)^3 (x^2 + 18x + 38)^4\\
& \ \ \times (x^2 + 32x + 65)^3 (x^2 + 82x + 9) (x^3 + 5x^2 + 90x + 31)\\
& \  (\textrm{modd} \ 113, J_{113}(x)^4).
\end{align*}
Comparing with the remainder of $\Phi_7(x^{113},x)$ modulo $J_{113}(x)^5$ shows that exactly the fourth power of $(x + 41) (x + 59) (x^2 + 18x + 38)$ divides $\Phi_7(x^{113},x)$, and this gives that
\begin{align*}
H_{-7 \cdot 113}(x) H_{-28 \cdot 113}(x) & \equiv x^4 (x+14)^4 (x+41)^8 (x+59)^8 (x^2 + 9x + 82)^6\\
& \times (x^2 + 18x + 38)^8 (x^2 + 32x + 65)^6 \ \ (\textrm{mod} \ 113).
\end{align*}
This agrees with the fact that $h(-7 \cdot 113) = h(-28 \cdot 113) = 32$.  This helps to explain the prevalence of the prime $113$ in the factorizations in Tables \ref{tab:1}-\ref{tab:3}.  Note that with $\sqrt{57} \equiv 31$ (mod $113$), we have the congruence
$$q_{171}(x) \equiv x^2 + 18x + 38 \ (\textrm{mod} \ 113),$$
so this gives an example where $q_{171}(x)$ divides $K_{7p}(x)$ modulo $p$.  Its conjugate factor is
$$\tilde q_{171}(x) \equiv (x + 14)(x + 41) \ (\textrm{mod} \ 113),$$
which does happen to divide $K_{7p}(x)$ in this case. \medskip

\noindent {\bf 5.} $p = 1217$.  Considering the gcd of $\Phi_7(x^{1217},x)$ and $J_{1217}(x)^3$ modulo $1217$ yields that
\begin{align*}
& H_{-7 \cdot 1217}(x) H_{-28 \cdot 1217}(x) \equiv (x+250)^2(x+941)^2 x^4 (x+477)^4 (x+765)^4\\
& \ \times (x+937)^4 (x+1168)^4\\
& \ \times (x^2 + 13x + 600)^4 (x^2 + 89x + 683)^4 (x^2 + 164x + 968)^4 (x^2 + 172x + 6)^4\\
& \ \times (x^2 + 263x + 931)^4 (x^2 + 413x + 546)^4 (x^2 + 463x + 537)^4 (x^2 + 515x + 390)^4\\
& \ \times  (x^2 + 786x + 754)^4 (x^2 + 805x + 1184)^4 (x^2 + 877x + 784)^4 \\
& \ \times  (x^2 + 1143x + 815)^4 (x^2 + 871x + 676)^6\\
& \ \times \big\{(x^2 + 14x + 190) (x^2 + 106x + 85) (x^2 + 257x + 897) (x^2 + 304x + 612)\\
& \ \times (x^2 + 307x + 276) (x^2 + 307x + 314) (x^2 + 410x + 94) (x^2 + 468x + 850)\\
& \ \times (x^2 + 478x + 1006) (x^2 + 522x + 299) (x^2 + 529x + 473) (x^2 + 535x + 576)\\
& \ \times (x^2 + 596 x + 566) (x^2 + 608x + 883) (x^2 + 656x + 307) (x^2 + 873x + 521)\\
& \ \times (x^2 + 944x + 560) (x^2 + 944x + 634) (x^2 + 1081x + 27) (x^2 + 1100x + 426)\\
& \ \times (x^2 + 1121x + 717)  (x^2 + 1129x + 1045)\big\}^2.
\end{align*}
This agrees with the class numbers $h(-7 \cdot 1217) = h(-28 \cdot 1217) = 110$.  From Table \ref{tab:3} we see that the occurrence of $(x^2 + 871x + 676)^6$ arises from the common factor of the polynomials
\begin{align*}
H_{-180}(x) & \equiv (x^2 + 786x + 754)(x^2 + 871x + 676),\\
H_{-192}(x) & \equiv (x^2 + 13x + 600)(x^2 + 871x + 676),\\
H_{-195}(x) & \equiv (x^2 + 877x + 784)(x^2 + 871x + 676)
\end{align*}
modulo $1217$.  The contributions of the remaining class equations, as well as the cofactors of this common factor, are correctly given in the factorization of Theorem 1.

\section{Appendix: Class equations and $Q_7(u,v)$.}

\begin{align*}
H_{-20}(x) &= x^2 - 1264000x - 681472000,\\
H_{-24}(x) &= x^2 - 4834944x + 14670139392,\\
H_{-40}(x) &= x^2 - 425692800x + 9103145472000,\\
H_{-48}(x) &= x^2 - 2835810000x + 6549518250000,\\
H_{-52}(x) &= x^2 - 6896880000x - 567663552000000,\\
H_{-75}(x) &= x^2 + 654403829760x + 5209253090426880,\\
H_{-115}(x) &= x^2 + 427864611225600x + 130231327260672000,\\
H_{-187}(x) &= x^2 + 4545336381788160000x - 3845689020776448000000,
\end{align*}

\begin{align*}
H_{-96}(x) &= x^4 - 23340144296736x^3 + 670421055192156288x^2\\
& \ \ + 447805364111967209472x - 984163224549635621646336,\\
H_{-132}(x) &= x^4 - 4736863498464000x^3 - 325211610485778048000000x^2\\
& \ \ + 54984539729717250048000000000x\\
& \ \ + 1656636925108948992000000000000,\\
H_{-160}(x) &= x^4 - 181195519824640800x^3 - 2940735389875294896000x^2\\
& \ \  + 13208221536779701382400000x - 293835053960432980416000000,\\
H_{-171}(x) &= x^4 + 694282057876537344x^3 + 472103267541360574464x^2\\
& \ \ + 8391550371275812148084736x - 1311901521779155773721411584,
\end{align*}
\begin{align*}
H_{-180}(x) &= x^4 - 2018504138609120000x^3 - 2867757758882006477169664000x^2\\
& \ \ + 16660473763558887652278272000000x\\
& \ \ - 6776421923145961044033929216000000,\\
H_{-192}(x) &= x^4 - 8041801037378436000x^3 + 15705521635909735050750000x^2\\
& \ \ + 826335556188178615474500000000x\\
& \ \ - 1080060886113159937649308593750000,\\
H_{-195}(x) &= x^4 + 11284411506057216000x^3 + 25349140792043819237376000x^2\\
& \ \ + 104773100319600336175104000000x\\
& \ \ - 233490285492432753672585216000000.
\end{align*}

\begin{align*}
& Q_7(u,v) = u^8 - 104545516658688000u^7 + (-34993297342013200v\\
& + 3643255017844740441130401792000000)u^6\\
 & + (-720168419610864v^2 - 1038063543615450389600613163008000v\\
 & - 42320664241971721884753245384947305283584000000000)u^5\\
 &  + (-4079701128594v^3 + 10685207605419643264415114978400v^2\\
 &  -  40689839325190046108805363111245786382336000000v\\
 & + 41375720005635744770247248526572116368162816000000000000)u^4\\
  & + (-9437674400v^4 - 16125487429364811901524078720v^3\\
 &  - 11269804822587811907034068119597607485440000v^2\\
  &  - 553081893983911854025093448165895580472573952000000000v\\
  &   - 13483958224762213714698012883865296529472356352000000000000000)u^3\\
  &    + (-10246068v^5 + 4460942463230217157721808v^4\\
   &    - 901645354876525685555163851940891808000v^3\\
   &     + 308881748688169626094151661068634506334830592000000v^2\\
     &    - 129852186866524354160683967940830111863092019200000000000000v\\
     &     + 1464765079488386840337633731737402825128271675392000000000000000000)u^2\\
     &      + (-5208v^6 - 177089350000162350352v^5 - 14066762315363594481667951670112000v^4\\
    &  - 17938541966217890687878494872256535904256000000v^3\\
     &  - 70610000800608793189892218223049143415629414400000000000v^2\\
    &    + 878989957472436106867912059683837492768381337600000000000000000v)u\\
     &    - v^7 + 312598951872417v^6 - 18309739022633325837478490400v^5\\
     &     + 89840545706297515985729629879602762400000v^4\\
       &    - 6015073802677018581437197939143986755770777600000000v^3\\
      &      + 212790112101924943054250805195683295626526720000000000000000v^2\\
       &      - 1708180850715319930422896480160236530761832857600000000000000000000v.
 \end{align*}

\medskip

\noindent Dept. of Mathematical Sciences, LD 270

\noindent Indiana University - Purdue University at Indianapolis (IUPUI)

\noindent Indianapolis, IN 46202

\noindent {\it e-mail: pmorton@iupui.edu}


\begin{thebibliography}{WWW}

\bibitem[1]{brm} J. Brillhart and P. Morton, Class numbers of quadratic fields, Hasse invariants of elliptic curves, and the supersingular polynomial, J. Number Theory 106 (2004), 79-111.

\bibitem[2]{co} D. A. Cox, {\it Galois Theory}, John Wiley \& Sons, 2004.

\bibitem[3]{f} R. Fricke, {\it Lehrbuch der Algebra, III: Algebraische Zahlen}, Friedr. Vieweg u. Sohn, Braunschweig, 1928.

\bibitem[4]{mor0} P. Morton, Explicit identities for invariants of elliptic curves, J. Number Theory 120 (2006), 234-271.

\bibitem[5]{mor1} P. Morton, Legendre polynomials and complex multiplication, I, J. Number Theory 130 (2010), 1718-1731.

\bibitem[6]{mor2} P. Morton, The cubic Fermat equation and complex multiplication on the Deuring normal form, Ramanujan J. 25 (2011), 247-275.

\bibitem[7]{mor5} P. Morton, Explicit congruences for class equations, Functiones Approx. 51 (2014), 77-110.

\bibitem[8]{mor4} P. Morton, The Hasse invariant of the Tate normal form $E_5$ and the class number of $\mathbb{Q}(\sqrt{-5l})$, J. Number Theory 227 (2021), 94-143.

\bibitem[9]{mor6} P. Morton, The Hasse invariant of the Tate normal form $E_7$ and the supersingular polynomial for the Fricke group $\Gamma_0^*(7)$, arXiv:2206.09801v2, 2022.

\end{thebibliography}
\end{document}